\documentclass[11pt]{amsart}

\usepackage[english]{babel}

\usepackage{color}
\usepackage{array}
\usepackage{enumerate}
\usepackage{enumitem}
\usepackage{refcount}
\usepackage{algorithm}
\usepackage{algorithmic}
\usepackage{boxedminipage}
\usepackage{float}
\usepackage[utf8x]{inputenc}
\usepackage{soul}

\usepackage{parskip}

\usepackage[colorlinks]{hyperref}
\hypersetup{
    colorlinks,
    linkcolor={red!50!black},
    citecolor={green!50!black},
    urlcolor={blue!80!black}
}
\usepackage[nameinlink,capitalize]{cleveref}
 
\usepackage{color}
\usepackage{graphicx}
\usepackage{multirow}
\usepackage{tikz}
\usepackage{amsmath}
\usepackage{amssymb}
\usepackage{caption}
\usepackage{tikz-cd}

\newtheorem{thm}{Theorem}[section]
\newtheorem{prop}[thm]{Proposition}
\newtheorem{lem}[thm]{Lemma}
\newtheorem{cor}[thm]{Corollary}

\theoremstyle{definition}
\newtheorem{defn}[thm]{Definition}
\newtheorem{nota}[thm]{Notation}

\newtheorem{remark}[thm]{Remark}

\def\proof{{\bf Proof.}}

\newcommand{\prfend}{\hbox to7pt{\hfil}

\par\vskip-\baselineskip\hbox to\hsize
{\hfil\vbox {\hrule width6pt height6pt}}\vskip\baselineskip}

\def\b{\par \noindent}

\def\dd{\medskip \par \noindent}
\long\def\eatit#1{}
\def\R{\mathbb{R}}

\def\N{\mathbb{N}}
\def\C{\mathbb{C}}

\def\P{\mathbb{P}}

\def\e{\varepsilon}

\def \"{``}
\font\tengothic=eufm10
\font\sevengothic=eufm7
\newfam\gothicfam
       \textfont\gothicfam=\tengothic
       \scriptfont\gothicfam=\sevengothic

\usepackage{xypic}

\DeclareMathOperator{\Res}{Res}
\DeclareMathOperator{\Tr}{Tr}

\begin{document}
\title{Postulation for 2-superfat points in the plane}

\author[S. Canino]{Stefano Canino}
\address[S. Canino]{Dipartimento di Matematica, Politenico di Torino, Corso Duca degli Abruzzi 24, 10129, Torino, Italy}
\email{stefano.canino@polito.it}

\author[M. V. Catalisano]{Maria V. Catalisano}
\address[M. V. Catalisano]{Scuola Politecnica, Universit\`a di Genova, Via All’Opera Pia 15, 16145, Genova, Italy}
\email{mariavirginia.catalisano@unige.it}

\author[A. Gimigliano]{Alessandro Gimigliano}
\author[M. Id\`a]{Monica Id\`a}
\address[A. Gimigliano, M. Id\`a]{Dipartimento di Matematica, Universit\`a di Bologna, Piazza di Porta S.Donato 5, 40126, Bologna, Italy}
\email{alessandr.gimigliano@unibo.it, monica.ida@unibo.it}

\author[A. Oneto]{Alessandro Oneto}
\address[A. Oneto]{Dipartimento di Matematica, Universit\`a di Trento, Via Sommarive, 14 - 38123 Povo (Trento), Italy}
\email{alessandro.oneto@unitn.it}

\maketitle

%%%%%%%%%%%%%%%%%%%%

\begin{abstract}
	We study the postulation of 0-dimensional schemes given by unions of $2$-superfat points in general position in the plane, i.e., the union of local schemes defined by the intersection of two distinct double lines. We prove that such schemes have good postulation, i.e., they have the expected Hilbert function. We also show the good postulation of such schemes when we add a general $3$-fat point. Finally, we use these results to answer a peculiar kind of interpolation problem.
\end{abstract}

%%%%%%%%%%%%%%%%%%%%

\section{Introduction}
In recent years the study of 0-dimensional schemes has seen a good deal of interest and also many applications (e.g. see \cite{BB},\cite{J},\cite{BBCG}) and there are still unanswered questions even in apparently simple cases (as the SHGH conjecture for the postulation of a set of generic fat points in the plane, e.g. see \cite{BCCGO}). The 0-dimensional schemes considered in this article sprang from ideas in a previous work of some of the authors. In \cite{CGI}, the study of the Jacobian ideal of a curve at an ordinary singular point lead us to consider possible structures of projective 0-dimensional {\it local} schemes, i.e., supported at one point, say $P$ in $\P^n$, which are {\it m-symmetric}, i.e., whose schematic intersection with any line through $P$ has the same length $m$. In \cite{CCGI}, the first four authors started a more accurate study of such $0$-dimensional schemes. In particular, we considered $m$-symmetric schemes that are maximal by inclusion, which we call {\it $m$-superfat points}, and we proved that their length is $m^n$.  In the case $m=n=2$, these schemes are defined by the schematic intersection of two double lines, i.e., their ideal is of the form $(L_1^2,L_2^2)$ where $L_1,L_2$ are two linear forms. We call these $2$-{\it squares}. 

	A first natural question about a given type $0$-dimensional projective scheme is about its {\it good postulation}, i.e., whether a generic union of them imposes independent conditions to forms of a given degree. The main result in the paper is that a general union of $2$-squares has good postulation, i.e., in algebraic terms, the defining ideal has maximal Hilbert function. See \Cref{Thetheorem}.

	\subsection*{Structure of the paper}
	In \Cref{sec:preliminaries}, we recall the basic definitions about superfat points and our main tools, such as the Horace method. In \Cref{sec:main_theorem}, we prove our main result \Cref{Thetheorem} on the good postulation of a general union of $2$-squares. In \Cref{sec:2squares_and_3fat}, we extend the result by adding a general $3$-fat point to the general union of $2$-squares. Finally, in \Cref{sec:interpolation}, we apply the above results to study a peculiar interpolation problem in the plane.

\section{Preliminaries}\label{sec:preliminaries}

\subsection{Notation}
	We work over the field of complex numbers $\C$. We denote by $S=\C[x,y,z]$ the coordinate ring of $\P^2$. We denote by $S_d$ the vector space of polynomials of $S$ having degree $d$, i.e., 
	\[
		S=\bigoplus_{d=0}^\infty S_d, \quad \dim S_d={d+2\choose 2}.
	\]
	Similarly, for any ideal $I \subset S$, we denote $I_d = S \cap I$. 
	
	Given a projective scheme $X \subset \P^2$, we denote by $I(X)$ the only homogeneous and saturated ideal of $S$ such that $X={\rm Proj}(S/I(X))$. If $X$ is $0$-dimensional, $\ell(X)$ denotes its length. The {\it Hilbert function} of $X$ is 
	\[
		H_X(d)=\dim_\C (S_d/I(X)_d).
	\]
	Recall that $I(X)_d$ corresponds to the vector space of curves of degree $d$ containing $X$.  If $X \subset r \subset \P^2$, where $r : \{L=0\}$ is a line, then we denote by $I_{X,r}\subset S/(L) $ the ideal of $X$ as a subscheme of $r$. 
	
	Finally, if $X_1,\dots,X_s$ are schemes of $\P^2$ with disjoint supports, we denote by $X_1+X_2+\dots+X_s$ their schematic union, while $mP$ denotes the {\it $m$-fat point} scheme associated to $I(P)^m$.

\subsection{Superfat points and the Horace method}\label{ssec:horace}

\begin{defn}
	A local 0-dimensional scheme $Z \subset \P^2$ supported at the point $P$ such that $\ell(Z \cap r) = m$ for any line $r$ passing through $P$ is called an {\bf $m$-symmetric scheme}. If $\ell(Z)$ is maximal, i.e. $\ell(Z) = m^2$, $Z$ is said to be an {\bf $m$-superfat point}.
\end{defn}

In \cite{CCGI}, the following explicit description of 2-superfat points in $\P^2$ is proved (while nothing similar exists for $m-$superfat points when $m\geq 3$). 

\begin{prop}
	Every 2-superfat scheme $Z\subset \P^2$ is such that $I(Z)=(x^2,y^2)$, up to projectivity.  We call such a scheme in $\P^2$ a {\bf $2$-square}.
\end{prop}
In the proof of \Cref{Thetheorem} we use {\it ``la méthode d’Horace''} that we briefly recall here and for which we refer to \cite{H}. Given a $0$-dimensional scheme $X \subset \P^2$ and a line $r:\{L=0\}$ we denote by 
\begin{itemize}
	\item $\Res_r(X)$: the \textbf{residue} of $X$ with respect to $r$, i.e. the scheme of $\P^2$ defined by the ideal $I(X) : (L)$;
	\item $\Tr_r(X)$: the \textbf{trace} of $X$ on $r$ given by the schematic intersection $X \cap r$, i.e., the subscheme of $r$ defined by the ideal $\frac{I(X)+(L)}{(L)}\subseteq S/(L)$.
\end{itemize}
The Horace method is based on the repeated use of the following exact sequence:
\begin{equation}\label{eq:exact_sequence}
	0 \rightarrow (I_{\Res_r(X)})_{d-1} \rightarrow (I_X)_d\rightarrow (I_{\Tr_r(X),r})_d .
\end{equation}
where $ (I_{\Tr_r(X),r})_d \subseteq S/(L)$  is $=\{0\}$ if only if $\ell(Tr_r(X),r) > d$, so (by Bezout) $r$ is a fixed component for the linear system of degree-$d$ curves through $X$, i.e. if and only if  $(I_{\Res_r(X)})_{d-1}\cong  (I_X)_d$. Therefore, if the goal is to show that $X$ has  the expected Hilbert function, i.e., by parameter count, $\max\{0, {d+2 \choose 2} - \ell(X)\}$ which is always a lower bound for $\dim (I_X)_d$, we have to find an iteration of specialisations of $X$ which, by passing to residue and trace, allows a successful induction on the degree $d$ which yields the required value for  $\dim (I_X)_d$. In order to do this, at every specialization and application of  \Cref{eq:exact_sequence} we must have not only that $\dim (I_{\Tr_r(X),r})_d=0$, but also that $\ell(Tr_r(X),r) = d+1 $, so to get that $\ell(\Res_r(X))= \ell(X) - (d+1)$ and $\dim(I_{\Res_r(X)})_{d-1}=\dim (I_X)_d$ has the required value. Once we get this results for a specialization of $X$ we will be able to conclude, since, by semicontinuity, the specialization of $X$  can only increase  $\dim (I_X)_d$.

%%%%%%%%%%%%%%%%%%%%%%%%%%%%%%%%%%%%%%%%  By semicontinuity, the specialization of the scheme can only increase the $h^0$.  $(I_{\Res_r(X)})_(d-1)\cong  (I_X)_d$

\section{Good postulation of general $2$-squares}\label{sec:main_theorem}
In this section we prove our main theorem on the Hilbert function of a general union of $2$-squares.

\begin{thm}\label{Thetheorem}
	If $X=Z_1+Z_2+\dots+Z_s\subseteq\P^2$ is a general union of $s$ 2-squares then
	\[
		H_X(d)=\min\left\{{d+2\choose 2}, \,4s\right\}
	\]
	or, equivalently, $ \dim I(X)_d = \max\left\{0, {d+2 \choose 2} - 4s\right\}$.
	% $$\dim I(X)_d={d+2\choose 2}-\min\{{d+2\choose 2}, \,4s\}=\max \{{d+2\choose 2}-4s, \, 0\} .$$
\end{thm}
\begin{remark}\label{rmk:extreme}
	A useful simplification in the proof of \Cref{Thetheorem} is that we can focus only on two extreme values of $s$:
	\begin{itemize}
		\item $s_*(d) := \left\lfloor \frac{{d+2 \choose 2}}{4}\right\rfloor$, i.e., the largest number of $2$-squares that we expect to impose independent conditions on the space of degree-$d$ plane curves. This because if a scheme $X$ imposes independent condition then the same holds for any $X' \subset X$.
		\item $s^*(d) := \left\lceil \frac{{d+2 \choose 2}}{4}\right\rceil$, i.e., the smallest number of $2$-squares that we expect to admit no degree-$d$ plane curve passing through them. This because if there are no curves through a scheme $X$ then the same holds for any $X' \supset X$.
	\end{itemize}
\end{remark}

\subsection{Our specializations} We introduce the specializations we will use in the proof of \Cref{Thetheorem}. First, we motivate our approach by underling the following arithmetic issue when trying to use the Horace method straightaway.

\begin{remark}\label{rmk:numerical_issues}
	As mentioned in \Cref{ssec:horace}, it is enough to prove the statement for a specialization of $X$. However, already for $d = 6$, we immediately encounter a numerical issue. Indeed, we have that $s^*(6) = s_*(6) = 7$ and the length of $X = Z_1 + \ldots + Z_7$ is exactly $28$ which matches exactly the dimension of plane sextics. We try to specialize $X$ to a scheme $\tilde{X}$ where some of the $2$-squares are supported along the same line, call it $r$. Our goal is to prove that $\dim(I_{\Res_r(\tilde{X})})_5 = \dim(I_{\Tr_r(\tilde{X}),r})_6 = 0$. In order to have $\dim(I_{\Tr_r(\tilde{X}),r})_6 = 0$ we need to specialize four $2$-squares to have support on $r$. However, in this way, the $\ell(\Res_r(\tilde{X})) = 20 < 21 = \dim(S_5)$, i.e., we would have $\dim(I_{\Res_r(\tilde{X})})_5 > 0$. The issue here is numerical and occurs whenever $d$ is even and $s_*(d) = s^*(d)$, i.e., whenever we cannot waste conditions on $\Tr_r(X)$ because otherwise we would be left with too little conditions on $\Res_r(X)$. 
\end{remark}

\begin{remark}\label{spec}
In light of \Cref{rmk:numerical_issues}, we need to specialise differently our scheme of $2$-squares when moving some of its components to have aligned supports on a line $r$. Without loss of generality, we suppose  $r: \{y=0\}$. Here we introduce a few types of specialisation.
\begin{enumerate}
	\item We can specialise a 2-square $Z$ in two ways:
	\begin{enumerate}
		\item  in order to have $$I(Z)=(x^2,y^2);$$
		\item in order to have $$I(Z)=((x-y)^2,(x+y)^2).$$
	\end{enumerate}
	In both cases, $\Res_y(Z)$ is a $0$-dimensional scheme of length $2$, but the difference is that: 
	\begin{itemize}
		\item in (a), $I(\Res_r(Z))=(x^2,y^2):(y)=(x^2,y)$. Hence, $\Res_r(Z)$ is a 2-jet contained in $r$ and $\deg(\Tr_r(\Res_r(Z)))=2$;
		\item in (b), $I(\Res_r(Z))=((x-y)^2,(x+y)^2):(y)=(x^2+y^2,xy,y^3):(y)=(x,y^2)$. Hence, $\Res_r(Z)$ is a 2-jet not contained in $r$ and $\deg(\Tr_r(\Res_r(Z)))=1$.
	\end{itemize}
	\item Given two $2$-squares with different support we may collapse them together. Namely, consider the family of pairs of $2$-squares depending on $t \in (0,1]$:
	\begin{align*}
	I(Z_t) & = (x^2,y^2) \cap ((x+tz)^2, (y+tz)^2) \\
	 & = \left(2x^{2}y-2xy^{2}+tx^{2}z-ty^{2}z, 2x^{3}-2y^{3}+3tx^{2}z-3ty^{2}z,y^{4}+2ty^{3}z+t^{2}y^{2}z^{2}\right).
	\end{align*}
	The limit scheme at $t = 0$ is a local $0$-dimensional scheme $Y$ of length $8$ defined by 
		\[
			I(Y) = \lim_{t \rightarrow 0} I(Z_t) = (x^2y-xy^2, x^3-y^3, y^4).
		\]
	Such scheme $Y$ has a ``slice" of length $3$ which allows us to overcome the arithmetic issue described in \Cref{rmk:numerical_issues} because allows us to have a trace of odd length.	 Let us describe the residues we get removing successively the line $r$. Indeed, we have that:
	\begin{itemize}
		\item $I(Y)+(y)=(x^3,y)$ so that $\deg(\Tr_r(Y))=3$ and \[I(\Res_r(Y))=I(Y):(y)=(x^2-xy,xy^2,y^3);\]
		\item $I(\Res_r(Y))+(y)=(x^2,y)$ so that $\deg(\Tr_r(\Res_r(Y)))=2$ and \[I(\Res_r(\Res_r(Y))) = I(\Res_r):(y)=(x,y)^2.\] Hence $\Res_r(\Res_r(Y))=2P$, where $P$ is the point $[0,0,1]$.
	\end{itemize}
\end{enumerate}
	Now let $m= \{x+\alpha y=0\}$ be a generic line through $P$. By an easy computation we get that, except for $\alpha = 0$ and $ \alpha = -1$, we have
	\begin{itemize}
	\item $I(\Res_r(Y))+(x+\alpha y)=(y^2,x+\alpha y)$, and $\deg(\Tr_{m}(\Res_r(Y)))=2$;
	\item  $\Res_m(\Res_r(Y))$ is defined by the ideal $I(\Res_r(Y)):(x+\alpha y)=(x,y)^2$, 
	that is, again, $\Res_m(\Res_r(Y))=2P$.
	\end{itemize}
	Notice that for $\alpha = 0$ or $ \alpha = -1$ we have $\deg(\Tr_{m}(\Res_r(Y)))=3$.
\end{remark}

\begin{nota}
	When in the following we use a specialisation of one of the types introduced in \Cref{spec} we refer to them respectively as specialisation of type (1.a), (1.b) or (2).
\end{nota}

\subsection{Lemmata}
Before proceeding with the proof of \Cref{Thetheorem}, we need some lemmata that will be used though our  inductive argument. 
\begin{lem}\label{dispari}
Let $d\in\N$ be odd and consider a $0$-dimensional scheme $\tilde X=X_1+X_2$ where:
\begin{itemize}
	\item $X_1=Z_1+Z_2+\dots+Z_{(d+1)/2}$, where all $Z_i$'s are $2$-squares with support on a line $r$ in such a way that $Z_1$ is specialised as type (1.b) and $Z_2,Z_3,\dots,Z_{{d+1}/2}$ are specialised as type (1.a);
	\item $X_2$ does not intersect $r$.
\end{itemize}
Then 
$$\dim I(\tilde X)_d=\dim I(X_2+P)_{d-2}$$
where $P$ is the support of $Z_1$.
\end{lem}
\proof ~ We apply the residual exact sequence in \Cref{eq:exact_sequence} to $\tilde{X}$. Since 
	\[
		\deg(\Tr_r(\tilde X))=2\cdot \frac{d+1}{2}=d+1,
	\]
	we have that $\dim(I_{\Tr_r(\tilde X),r})_d = d+1 - (d+1) = 0$; hence, $\dim(I_{\tilde X})_d = \dim(I_{\Res_r(\tilde X)})_{d-1}$. Now, we are left with 
	\[
		X':= \Res_r(\tilde X) = \Res_r(X_1)+\Res_r(X_2)=Z_1'+Z_2'+\dots+Z_{(d+1)/2}'+X_2,
	\]
	with $Z_i':=\Res_r(Z_i)$ for $i=1,\dots,\frac{d+1}{2}$. Note that $\Res_r(X') = X_2 + P$. Now, as explained in \Cref{spec}, we have that
	\[
		\deg(\Tr_r(X'))=\sum_{i=1}^{(d+1)/2}\deg(\Tr_r(Z_i'))=1+2 \left(\frac{d+1}{2}-1\right)=d.
	\]
	Thus, $\dim(I_{Tr_r(X'),r})_{d-1} = d - d = 0$ and, by the residual exact sequence in \Cref{eq:exact_sequence} applied to $X'$, we get that $\dim(I_{\Res_r(\tilde X)})_{d-1}=\dim(I_{X_2+P})_{d-2}.$ This concludes the proof.
\prfend 

\begin{lem}\label{24}
	Let $X=Z_1+Z_2+\dots+Z_s$ be a general union of $s$ $2$-squares in $\P^2$ and let $d\in\{2,4\}$. Then
	\[
		\dim I(X)_d=\max\left\lbrace 0,{d+2\choose 2}-4s\right\rbrace.
	\]
\end{lem}
\proof ~ As recalled in \Cref{rmk:extreme}, for each $d$, it is enough to consider $s \in \{s^*(d),s_*(d)\}$. Hence, it is enough to analyze the following four cases.
\begin{itemize}[leftmargin=*]
	\item $d=2$ and $s=1$\dd
		This is immediate because a single $2$-square is not contained in a line, i.e., $\dim I(X)_1=3$, and it has length $4$, i.e., $\dim I(X)_d=4$ for any $d\geq 2$.

	\item $d=2$ and $s=2$\dd
		As explained in \Cref{spec}, it is enough to prove the statement for a specialization of $X$. We consider the $0$-dimensional scheme obtained by collapsing the two $2$-squares $Z_1$ and $Z_2$ as in specialization (2). Since the ideal of such limit does not have generators of degree $2$, we are done.

	\item $d=4$ and $s=3$\dd
		We assume that $X$ is specialized as $X = Y + Z_3$, where $Y$ and $Z_3$ have support on the line $r$, $Y$ is obtained by collapsing the two 2-squares as in specialization (2) and $Z_3$ is specialized as in specialization (1.a) of \Cref{spec}. Since $\deg(\Tr_r(X))=3+2=5$, then $\dim(I_{\Tr_r(X),r})_4= 0$ and, by \Cref{eq:exact_sequence}, $\dim(I_X)_4=\dim(I_{\Res_r(X)})_3$. We are left with \[X':=\Res_r(X)=Y'+Z_3'\] where $Y':=\Res_r(Y)$ and $Z_3':=\Res_r(Z_3)$. Note that \[\deg(\Tr_r(X'))=\deg(\Tr_r(Y'))+\deg(\Tr_r(Z_3'))=2+2=4\] and $\Res_r(X')=\Res_r(Y')=2P$, where $2P$ is a $2$-fat point with support at $P \in r$,. Therefore, $(I_{\Tr_r(X'),r})_3=\{0\}$ and, by \Cref{eq:exact_sequence} applied to $X'$, we get that $\dim(I_{X'})_3=\dim(I_{2P})_2= 3$ which matches the expected dimension ${4 + 2 \choose 2} - 3\cdot 4$.

	\item $d=4$ and $s=4$\dd 
		Let $X = Y + Z_3 + Z_4$, where $Y$ and $Z_3$ are as in the latter case and $Z_4$ is a general $2$-square. Proceeding exactly as in the latter case, we find that $\dim(I_X)_4=\dim(I_{2P+Z_4})_2$. Since $\ell(\Tr_r(2P+Z_4)) = 4$, then $\dim(I_{\Tr_r(2P+Z_4)})_2 = 0$ and, by \Cref{eq:exact_sequence} applied to the scheme $2P+Z_4$, we further reduce to $\dim(I_{2P+Z_4})_2 = \dim(I_{\Res_r(2P+Z_4)})_1$ with $\Res_r(2P+Z_4) = P + Z'_4$ and $\ell(Z'_4) = 2$. Since $Z_4$ is assumed to be general, $Z'_4$ is not contained in the line $r$. Hence, we conclude that $(I_X)_4 = (I_{P+Z'_4})_1 = 0$ as expected. \prfend
\end{itemize} 

\begin{lem}\label{pari}  
	Let $d$ be an even integer, $d\geq 6$, and let $r$ and $m$ be two general lines, with $P=r\cap m$. Let $Z_1,\dots, Z_d$ be general $2$-squares and consider the scheme  $\tilde X=X_1+X_2+X_3$ where:
\begin{itemize}
	\item $X_1=Y+Z_3+\dots+Z_{d/2 +1}$, where the scheme $Y$ has support on $P$ and it is obtained by collapsing together  $Z_1$ and $Z_2$ as in specialisation (2), and the $Z_i$'s, $3 \leq i \leq d/2 +1$,  have support on  $r\smallsetminus \{P\}$,
$Z_3, \dots,Z_{d/2}$ are specialised on $r$ as type (1.a), and $Z_{d/2 +1}$ is specialised as type (1.b) ;
	\item $X_2=Z_{d/2 +2}+\dots+Z_{d}$, where the $Z_i$'s have support on  $m\smallsetminus \{P\}$, $Z_{d/2 +2},\dots,Z_{d-1}$ are specialised on $m$ as type (1.a), and $Z_{d}$ is specialised as type (1.b);
	\item  $X_3$ does not intersect $r\cup m$.
\end{itemize}
Then
	\[
		\dim I(\tilde X)_d=\dim I(P_1+P_2+X_3)_{d-4}
	\] 
	where $P_1$ and $P_2$ are the supports of $Z_{d/2+1}$ and $Z_d$, respectively.
\end{lem}
\proof ~ By construction,
	\[
		\deg(\Tr_r(\tilde X))=2\left(\frac{d}{2}-1\right)+3=d+1.
	\]
	Thus, $\dim(I_{\Tr_r(\tilde X),r})_d=0$ and, by \Cref{eq:exact_sequence} applied to $\tilde{X}$ with respect to the line $r$, we get $\dim(I_{\tilde X})_d=\dim(I_{\Res_r(\tilde X)})_{d-1}$. We are left with
	\[
		X':=\Res_r(\tilde X)=\Res_r(Y)+\Res_r(Z_3+\dots+Z_{d/2+1})+X_2+X_3
	\]
	Since
	\[
		\deg(\Tr_m(X'))=2+2 \left(\frac{d}{2}-1\right)=d
	\]
	we get $(I_{\Tr_m(X'),m})_{d-1}=\{0\}$ and, by applying \Cref{eq:exact_sequence} to $X'$ with respect to the line $m$, we deduce $\dim(I_{X'})_{d-1}=\dim(I_{\Res_m(X')})_{d-2}$, where 
	\[
		X'':=\Res_m(X')=\Res_m(\Res_r(Y))+\Res_r(Z_3+\dots+Z_{d/2 +1})+\Res_m(Z_{d/2 +2}+\dots+Z_{d})+X_3.
	\]
	By \Cref{spec}, $\Res_m(\Res_r(Y)) =2P$ so we get
	\[
		\deg(\Tr_r(X''))=2+2\left (\frac{d}{2}-2 \right )+1=d-1
	\]
	and $\dim(I_{\Tr_r(X''),r})_{d-2}=0$. Then, by applying \Cref{eq:exact_sequence} to $X''$ with respect to the line $r$, we deduce $\dim(I_{X''})_{d-2}=\dim(I_{\Res_r(X'')})_{d-3}$, where 
	\begin{align*}
		X'''& :=\Res_r(X'') \\
		& =\Res_r(2P))+\Res_r(\Res_r(Z_3+\dots+Z_{d/2 +1}))+\Res_m(Z_{d/2 +2}+\dots+Z_{d})+X_3 \\
		& = P + P_1+\Res_m(Z_{d/2 +2}+\dots+Z_{d})+X_3.
	\end{align*}
	Now, 
	\[
		\deg(\Tr_m(X'''))=1+2\left (\frac{d}{2}-2 \right )+1=d-2
	\]
thus $I_{\Tr_m(X'''),m})_{d-3}=\{0\}$ and, by \Cref{eq:exact_sequence} applied to $X'''$ with respect to the line $m$, we get $\dim(I_{X'''})_{d-3}=\dim(I_{\Res_m(X''')})_{d-4}$, where $\Res_m(X''')= P_1+P_2+X_3$.
\prfend

\subsection{Proof of \Cref{Thetheorem}} We have to prove that, 
if $X=Z_1+Z_2+\dots+Z_s$ is a general union of $s$ 2-squares  in $\P^2$, then
$$\dim I(X)_d=\max\left\lbrace{d+2\choose 2}-4s, \, 0\right\rbrace.$$

 By \Cref{rmk:extreme}, for each degree $d$, we consider $s \in \{s_*(d), s^*(d)\}$. The proof is by induction on $d$: we prove that if the statement holds for $d-2$ when $d$ is odd and for $d-4$ when $d$ is even, then the thesis is true for $d$.

 The initial cases are $d=1$ (trivial), $d=2$ and $d=4$ (see \Cref{24}).  

 Notice that, for any $d\in\N$,  there are $\e',\, \e'' \in\{0,1,2,3\}$ such that 
$$s_*(d)=\left\lfloor {1\over 4}{d+2 \choose 2}\right\rfloor={1\over 4}{d+2 \choose 2}-{\e'\over 4}, \qquad \quad s^*(d)=\left\lceil {1\over 4}{d+2 \choose 2}\right\rceil={1\over 4}{d+2 \choose 2}+{\e''\over 4}.$$
We distinguish four cases according to the parity of $d$ and $s=s_*(d)$ or $s=s^*(d)$.

\begin{itemize}[leftmargin=*]
\item Case 1: $d \geq 3$ odd and $s=s_*(d)={1\over 4}{d+2 \choose 2}-{\e \over 4}$.

	We want to prove that $\dim I(X)_d=\max\left\lbrace{d+2\choose 2}-4s, \, 0\right\rbrace=\e$. We specialise $X$ to a scheme $\tilde X = X_1 + X_2$ as in the statement of \Cref{dispari} with $X_2$ generic union of 
	$$s-{d+1\over 2}= {1\over 4}{d+2 \choose 2}-{\e \over 4}-{d+1\over 2}= {1\over 4}\left({d \choose 2}-1-\e \right)$$ 2-squares. \Cref{dispari} tells us that
	$\dim I(\tilde X)_d=\dim I(X_2+P)_{d-2}$,  where $P$ is a generic point. By induction assumption we have 
	$$\dim I(X_2)_{d-2}=\max\left\lbrace{d\choose 2}-4{1\over 4}\left({d \choose 2}-1-\e\right), \, 0\right\rbrace=1+\e,$$
	hence $\dim I(X_2+P)_{d-2}=\e$, so we get $\dim I(X)_d\leq \dim I(\tilde{X})_d=\e$. On the other hand, the expected Hilbert function is a lower bound; that is, $\dim I(X)_d\geq {d+2 \choose 2}- 4\left({1\over 4}{d+2 \choose 2}-{\e \over 4}\right)= \e$.

\item Case 2: $d\geq 3$ odd and $s=s^*(d)={1\over 4}{d+2 \choose 2}+{\e \over 4}$.

	If $\e = 0$, $s_*(d)=s^*(d)$ and we are in Case 1. Hence, we assume $\e>0$. We have to prove that 
	$\dim I(X)_d = 0$.
	We specialise $X$ to a scheme $\tilde{X}$ similarly to Case 1, except that now $X_2$ is the generic union of 
	$$s-{d+1\over 2}= {1\over 4}{d+2 \choose 2}+{\e \over 4}-{d+1\over 2}= {1\over 4}\left({d \choose 2}-1+\e \right)$$ 2-squares.
	By induction assumption we have 
	$$\dim I(X_2)_{d-2}=\max\left\lbrace{d\choose 2}-4{1\over 4}\left({d \choose 2}-1+\e\right), \, 0\right\rbrace=0.$$
	Hence, $\dim I(X)_d\leq \dim I(\tilde{X})_d=\dim I(X_2+P)_{d-2}=0$, where $P$ is any point, and we are done.

\item Case 3: $d\geq 6$ even and $s=s_*(d)= {1\over 4}{d+2 \choose 2}-{\e\over 4}$.

	We have to prove that $\dim I(X)_d=\e.$ We specialise $X$ to a scheme $\tilde X$, where $\tilde{X}$ is as in the statement of \Cref{pari}, with $X_3$ general union of 
	$$s-d={1\over 4}\left({d+2 \choose 2}-\e \right)-d= {1\over 4}\left({d-2 \choose 2}-2-\e \right)$$ 
	2-squares. \Cref{pari} tells us that
	$\dim I(\tilde X)_d=\dim I(X_3+P_1+P_2)_{d-4}$, where $P_1$ and $P_2$ are general points. By induction assumption we have 
	$$\dim I(X_3)_{d-4}=\max\left\lbrace{d-2\choose 2}-4{1\over 4}\left({d-2 \choose 2}-2-\e \right), \, 0\right\rbrace=2+\e,$$
	hence $\dim I(X_3+P_1+P_2)_{d-4}=\e$ and $\dim I(X)_d\leq \dim I(\tilde{X})_d=\e$. This matches the lower bound given by the expected Hilbert function and it concludes the proof. 

\item Case 4: $d\geq 6$ even and $s=s_*(d)={1\over 4}{d+2 \choose 2}+{\e \over 4}$.

	We have to prove that $\dim I(X)_d=0.$ We specialise $X$ to $\tilde{X}$ as in Case 3, except that now $X_3$ is the generic union of 
	$$s-d= {1\over 4}\left({d-2 \choose 2}-2+\e \right)$$ 2-squares.
	By induction assumption we have 
	$$\dim I(X_3)_{d-4}=\max\left\lbrace{d-2\choose 2}-4\cdot {1\over 4}\left({d-2 \choose 2}-2+\e \right), \, 0\right\rbrace=\max\{\e-2,0\}$$
	and thus, since $0\leq \e\leq 3$, we get $\dim I(X_3+P_1+P_2)_{d-4}=0$. Since, by \Cref{pari}, we have $\dim I(X)_d\leq \dim I(\tilde{X})_d=\dim I(X_3+P_1+P_2)_{d-4}=0$, we conclude the proof. \prfend
\end{itemize}

\begin{remark}\label{rmk:marvi}
In  case  $d$ even, and $s= s_*(d)=s^*(d) $, we have found also another way to  overcome the arithmetic issue described in \Cref{rmk:numerical_issues}, and we would like to show to the reader a sketch of this alternative method, which avoids the use of collisions.

 The proof of Theorem \ref{Thetheorem} is again by induction on $d$, but when we deal  with the problematic case,  instead of considering the scheme $X$, we first study a scheme $Y$ union of $s-1$ $2$-squares, and a double point $2Q$.  Then 
we consider two general lines $r$ and $m$, and  we specialize $Y$ to 
 the following scheme
$\tilde Y=X_1+X_2+X_3$, where  

\begin{itemize}
      
	\item $X_1=Z_1+Z_2+\dots+Z_{d/2}$, where all $Z_i$'s are $2$-squares with support on the line $r$ in such a way that $Z_1$ is specialised as type (1.b) and $Z_2,Z_3,\dots,Z_{d/2}$ are specialised as type (1.a);
	\item $X_2=W_1+W_2+\dots+W_{d/2}$, where all $W_i$'s are $2$-squares with support on the line $m$ in such a way that $W_1$ and  $W_2$ are specialised as type (1.b) and $W_3,\dots,Z_{{d}/2}$ are specialised as type (1.a);
	\item $X_3$ is union of $(s-d-1)$ general $2$-squares and the double point $2Q$.
\end{itemize}
Now let $P=r \cap m$. By Bezout's Theorem we  may remove the lines $r$ and $m$ two times from the curves of the linear system $I(\tilde Y+P)_d $, and so $\dim I(\tilde Y+P)_d =\dim I(X_3+ P_1+P_2+P_3)_{d -4}$, where $P_1, P_2, P_3$ are  generic points (they are the supports of $Z_1, W_1, W_2$, respectively).  

Since, by the induction assumption,  $(s-d)$ $2$-squares  give the right number of conditions to the curves of degree $d-4$, hence the same happens for $(s-d-1)$ $2$-squares and a 2-fat point, we get
\[
\dim I(\tilde Y+P)_d =
\dim I(X_3+ P_1+P_2+P_3)_{d-4}= {d-2\choose 2} -4(s-d-1)-3-3=0,
\]
and so also $\dim I(Y+P)_d =0$.
It follows that  $\dim I(Y )_d =1$, that is, there exists  only one curve  $C$  of degree $d$ containing $(s-1)$ $2$-squares, and a double point $2Q$.

If there exists a 2-square with support on $Q$, not contained in $C$, then
going back to $X$, the dimension  will decrease by one, and we are done. So assume that $C$ contains every  2-square with support on $Q$.
By \cite{CCGI} we get that $C$ has a triple point in $Q$. Hence by simmetry, the same shoud happen at every other 2-square of the scheme, i.e. $C$ should have $s$ triple points instead of $s$ 2-squares, which is impossible (e.g. see \cite{M}, \cite{GI}).

 \end{remark}

%%%%%%%%%%%%%%%%%%%%%%%%%%%%%%%%%%%%%%%%

\section{Postulation of a general set of $2$-squares and a $3$-fat point}\label{sec:2squares_and_3fat}
	\begin{nota}In this section, we denote by
	\[
		X_s:=Z_1+\dots+Z_s+3Q
	\]
	the generic union of $s$ 2-squares $Z_1,\dots,Z_s$ and a triple point $3Q$. Let $P_i$ be the support of $Z_i$. 
	\end{nota}
	This section is devoted to prove the following result.
	\begin{thm}\label{2superfats+triple} 
		For any $(s,d)\neq(1,3)$, the dimension of $I(X_s)_d$ is the expected one, i.e.
		$$\dim I(X_s)_d = \max\left\{{d+2 \choose 2} - 4s - 6, 0\right\}.$$
		In the case $X_1=Z_1+3Q$ we have $ \dim I(X_2)_3 = 1$ and not $0$ as expected.
	\end{thm}
	The proof of \Cref{2superfats+triple} will proceed by induction on $d$. The base cases are collected in \Cref{initialcases}. The induction step is different depending on the parity of $d$: in \Cref{paritriplo} we consider the case $d$ even, while in \Cref{disparitriplo} we consider the case $d$ odd. 

	\begin{nota} 
		We set the following statement:

		$A(s,d):$ \lq\lq The dimension of $I(X_s)_d$ is the expected one, i.e. $\dim(I(X_s)_d) = \max\{\dim S_d - 4s - 6, 0\}$\rq\rq.
	\end{nota}

	As motivated in \Cref{rmk:extreme}, we will focus on the following two values of $s$:
	\[
		s_*(d)=\left\lfloor {{d+2 \choose 2}-6\over 4}\right\rfloor \quad\text{and}\quad s^*(d)=\left\lceil {{d+2 \choose 2}-6\over 4}\right\rceil.
	\]
% which are, respectively, the maximum $s$ such that $\dim(I_{X_s})_d>0$ and the minimum $s$ such that $\dim(I_{X_s})_d\leq0$. 
\begin{lem}\label{initialcases} 
	$A(s,d)$ is true for any $(s,d)$ with $1\leq d\leq 5$, except for $(s,d)=(1,3)$ where we have $\dim I(X_1)_3=1$.
\end{lem}
 \proof ~ The cases $d\in\{1,2\}$ are trivial because $3Q$ is not contained in any line or conic. We distinguish now the remaining cases.
\begin{itemize}[leftmargin=*]
	\item If $d=3$ and $s=1$, it is immediate to see that we have the exceptional case for which $I(X)_3 = \langle L^3 \rangle$ where $L$ is the line passing through $P_1$ and $Q$. For $s \geq 2$, since, by genericity assumption, $P_2$ does not lie on the line $\overline{P_1Q}$, we have that $\dim(I(X_s)_3)=0$.

	\item For $d=4$, we consider only $s_*(4)=2$ and $s^*(4)=3$. If $s=2$, then the lines $r_1=\overline{P_1Q}$ and   $r_2=\overline{P_2Q}$ are fixed components for $I(X_3)_4$ because both intersect $X_3$ with multiplicity $5$. Then,  $I(X_2)_4 = (L_1L_2) \cdot I(X'_3)_2$, where $r_i=\{L_i=0\}$ and $X'_3$ is the union of $Q$ and two general $2$-jets, i.e., schemes of length $2$, supported at $P_1$ and $P_2$. It is well-known that $X'_3$ has generic Hilbert function, i.e. there is an unique conic $\{C=0\}$ passing through $P_1,P_2,Q$ with assigned tangents at $P_1$ and $P_2$.  Hence, $\dim I(X_3)_4 = \dim I(X'_3)_2 = 1$, and the only form in $I(X_3)_4$ is $L_1L_2C$, as expected. Moreover, this implies that for $s=3$ we have $\dim I(X_3)_4 = 0$ by genericity assumption.

	\item For $d=5$ we consider the cases $s_*(5)=3$ and $s^*(5)=4$. If $s=3$ we fix a line $r$ and we specialise $X_3$ to $X_3'=Z_1'+Z_2'+Z_3'+3Q$, where $Z_1'$ and $Z_2'$ are supported $r$ as type (1.a) of \Cref{spec}, $Z_3$ is supported on $r$ as type (1.b) and $Q$ does not intersect $r$. Using twice the residual exact sequence \Cref{eq:exact_sequence} with respect to $r$, we find that $\dim I(X_3')_5=\dim I(X_3'')_3 =3$ where $X_3''=\Res_r(\Res_r(X_3'))=P_3'+3Q$. Thus, we conclude by semicontinuity. If $s=4$, by using an analogous specialisation $X_4'$ we find $\dim I(X_4')_5=\dim I(X_4'')_3$ with $X_4'''=P_3'+Z_4+3Q$. Now, from the case $d=3$, we know that there is only one cubic containing $Z_4$ and $3Q$. By genericity of $P_3'$, we deduce that $\dim I(X_4')_5=0$ and the result follows again by semicontinuity. \prfend
\end{itemize}  

\begin{prop}\label{paritriplo}
	$A(s,d)$ is true for any $(s,d)$ with $d$ even and $d\geq 6$.
\end{prop}
\proof ~ Let us set $d=2k$, with $k\geq 3$. We have to prove that 
	\[
		\dim I(X_s)_d=\max\left\{{d+2 \choose 2}-4s-6,0\right\}=\max\{2k^2+3k-5-4s,0\}.
	\]
	Note that a straightforward computation shows that $s^*(d)\geq s_*(d)\geq k-1$ for any $k\geq 3$.
	%  and that  {\color{red} allows} us to treat the cases $s=s_*(d)$ and $s=s^*(d)$ simultaneously. 
	This allows us to consider the same specialization for both $s \in \{s_*(d),s^*(d)\}$, namely
	\[
		X_s' = Z_1'+\dots+Z_{k-1}'+Z_k+\dots Z_s+3Q
	\]
	where $Z_1',\dots,Z_{k-1}'$ are specialised on a line $r$ with specialisation of type (1.a) as in \Cref{spec}, $Q\in r$ and $Z_k,\dots,Z_s$ do not meet the line $r$. We have 
	$$\deg(\Tr_r(X_s'))=2(k-1)+3=2k+1=d+1$$ 
	and
	$$\deg(\Tr_r(\Res_r(X_s')))=2(k-1)+2=2k=d.$$ 
	Thus, using twice the residual exact sequence as in \Cref{eq:exact_sequence} with respect to $r$, we find
	$$\dim(I(X_s'))_d=\dim(I(\Res_r(\Res_r(X_s'))))_{d-2}$$
	where $\Res_r(\Res_r(X_s'))=Z_k+\dots Z_s+Q=:X_s''$. Since $X_s''$ is a general union of $(s-k+1)$ 2-squares and a simple point, then $I(X_s'')_{d-2}$ has, by \Cref{Thetheorem} the expected dimension, i.e.
	$$\dim I(X_s'')_{d-2}=\max\left\{{d \choose 2}-4(s-k+1)-1,0\right\}=\max\{2k^2+3k-5-4s,0\}.$$
	The statement follows by semicontinuity. \prfend

\begin{prop}\label{disparitriplo}
	$A(s,d)$ is true for any $(s,d)$ with $d$ odd and $d\geq 7$.
\end{prop}
	\proof ~ Let us set $d=2k-1$ with $k\geq 3$. We have to prove that
	$$\dim I(X_s)_d=\max\left\{{d+2 \choose 2}-4s-6,0\right\}=\max\{2k^2+k-6-4s,0\}.$$
	We prove the result by induction on $d$, i.e. we suppose that $A(s,d-2)$ is true for any $s$ and we show that $A(s-k,d-2)$ implies $A(s,d)$. We use as base case $A(s,5)$ which is true by \Cref{initialcases}.
	It is immediate to see that $s^*(d)\geq s_*(d)\geq k$; hence, for any $s \in \{s_*(d), s^*(d)\}$, we specialise $X_s$ to
	$$X_s'=Z_1'+\dots+Z_{k}'+Z_{k+1}+\dots Z_s+3Q$$
	where $Z_1',\dots,Z_{k-1}'$ are specialised on a line $r$ with specialisation (1.a) of \Cref{spec}, $Z_k'$ is specialised on $r$ with specialisation (1.b) and $Z_{k+1},\dots,Z_s, 3Q$ are generic and away from $r$. We have
	$$\deg(\Tr_r(X_s'))=2k=d+1$$ 
	and
	$$\deg(\Tr_r(\Res_r(X_s')))=2(k-1)+1=2k-1=d.$$ 
	Thus, using twice the residual exact sequence in \Cref{eq:exact_sequence} with respect to $r$, we find
	$$\dim(I(X_s'))_d=\dim(I(\Res_r(\Res_r(X_s'))))_{d-2}$$
	where $\Res_r(\Res_r(X_s'))=P_k'+Z_{k+1}+\dots Z_s+3Q=:X_s''$. By induction $X_s''$ has good postulation in degree $d-2$ and thus
	$$\dim I(X_s'')_{d-2}=\max\left\{{d \choose 2}-4(s-k)-6-1,0\right\}=\max\{2k^2+k-6-4s,0\}.$$
	The statement follows by semicontinuity. \prfend

%%%%%%%%%%%%%%%%%%%%%%%%%%%%%%%%%%%%%%%%

\section{An interpolation problem}\label{sec:interpolation}

	We interpret here our previous results in terms of the following interpolation problem.
	
	Given $s$ general points $P_1,...,P_s$ in the plane and a {\it frame of reference} at each of them,  
	i.e. two distinct lines $\{r_{i1},r_{i2}\}$ through each $P_i$, we want to study the algebraic curves of degree $d$ having a singular double point at each $P_i$ which is either a node,  and its two tangents are symmetric with respect to $\{r_{i1},r_{i2}\}$ (i.e. their equations can be written as $a_iL_{i1}+b_iL_{i2}=0$ and $a_iL_{i1}-b_iL_{i2}=0$), where $r_{ji}=\{L_{ji}=0\}$, while, if it has a double tangent line, it is either $r_{i1}$ or $r_{i2}$. 

	The given conditions are associated to a linear system $V_d$ of curves in the plane. We want to know what its dimension is and whether its generic curve satisfies the required hypothesis (it could be more singular). In order to answer these questions we consider the problem in the projective plane, and show that the required conditions are associated to the postulation of a generic union of $2$-squares.

\begin{prop}\label{solodoppi} 
	Let $V_d=I(Z_1 +\ldots + Z_{s-1} + Z_s)_d$, where the $Z_i$'s are generic $2$-squares with $I_{Z_i} = (L_{i1}^2,L_{i2}^2)$. If the dimension of $V_d$ is positive, then its generic element is a curve which has only double points at each $P_i$ that are either ordinary nodes whose tangents are symmetric with respect to $L_{i1},L_{i2}$, or $\dim V_d=1$ and the only curve in it has $L_{ij}^2$ as double tangent cone at some $P_i$ (maybe at all of them).
\end{prop}
 \proof ~ \Cref{2superfats+triple} implies that when we consider a linear system $V_d = I(Z_1+\ldots + Z_s)_d$, its generic element gives a curve that has only double points at each $Z_i$. In fact if we were to have that the generic element in $V_d$ has multiplicity at least 3 at one point, say $P_s$, then it would belong also to $I(Z_1+\ldots +Z_{s-1}+3P_s)_d$, and this cannot happen since, by \Cref{2superfats+triple}, $I(Z_1+\ldots +Z_{s-1}+3P_s)_d$ has dimension $\max\{0, \dim I(Z_1+\ldots+Z_s)_d-2\}$. 

 Now, since at each $P_i$ we have $I_{Z_i} = (L_{i1}^2,L_{i2}^2)$, without loss of generality, we suppose $I_{Z_1} = (x^2,y^2)$. We have seen that a generic element of $V_d$ has actually multiplicity $2$ at $P_1$, hence its initial part is of type $a^2x^2-b^2y^2 = (ax+by)(ax-by)$, $a,b\in \C$. Thus, either $ab\neq 0$, and the curve has tangent cone of type $r_ir'_i$ in $P_i$, with the two tangent lines which are symmetric with respect to the axes, or $ab=0$ and the tangent cone is either $x^2$ or $y^2$, as required. Similarly, for all other $Z_i$'s. 
 
 Moreover, this cannot happen if $\dim V_d \geq 2$. Actually, in this case, if we had that the generic element in $V_d$ has its initial part at some $P_i$, say for $i=s$, of type $x^2 + ({\it high~order~terms})$, it would follow that a basis of $V_d$ can be written as $\{f_1,...,f_s\}$, $s\geq 2$, with $f_1 = x^2 + ({\it high~order~terms})$ while the $f_i$'s for $i\geq 2$ have initial degree $\geq 3$, but in this case if we consider the vector space $V'_d = \langle f_2,...,f_s\rangle$, we have $\dim V'_d = \dim V_d-1$, but this is impossible since $ V'_d \subset 
  I(Z_1+\dots+Z_s-1,3P_s)_d$ and $\dim I(Z_1+\dots+Z_{s-1},3P_s)_d = \dim V_d-2$ by \Cref{Thetheorem}.
\prfend

Let us also notice that if $\dim V_d \geq 2$ and we impose to the elements of $V_d$ that its tangent cone at a point $P_i$ contains a fixed generic line $r_i$ through $P_i$, then the dimension of the system drops exactly by 1 (otherwise, it could not drop by $2$ when we impose $3P_i$). Hence, the curves with tangent cone $r_ir'_i$, where $r'_i$ is the  line symmetric to $r_i$ with respect to $L_1,L_2$, form a dense open set in a hyperplane of $V_d$, which has to intersect the subvector $\R$-space given by the equations with real coefficients. So, if we consider the whole situation over the field $\R$, the linear systems we considered have the same dimensions, and the curve with given tangent cone $r_ir_i'$ does exist also in the real case (while if $\dim V_d =1$ the only curve in it could have isolated singularities at some $P_i$). In summary, we obtain the following corollary.

\begin{cor} 
	If $\dim V_d\geq 2$, \Cref{solodoppi} holds also over the field $\mathbb{R}$.
\end{cor}  

\medskip \b {\bf Funding.} AG proposed the problem to AO during the semester AGATES in Warsaw where they were partially supported by the Thematic Research Programme Tensors: geometry, complexity and quantum entanglement, University of Warsaw, Excellence Initiative - Research University and the Simons Foundation Award No. 663281 granted to the Institute of Mathematics of the Polish Academy of Sciences for the years 2021-2023.

SC  has been partially funded by the European Union under NextGenerationEU, Mission 4 Component 2, 2022E2Z4AK, PRIN 2022. AG and MI have been partially funded by the EuropeanUnion under NextGenerationEU, Mission 4 Component 2, J53D23003750006, PRIN 2022 - Prot. n. 2022E2Z4AK -  0-Dimensional Schemes, Tensor Theory, and Applications. AO has been partially funded by the European Union under NextGenerationEU, Mission 4 Component 2, E53D23005400001, PRIN 2022 - Prot. n.20223B5S8L - Birational geometry of moduli spaces and special varieties.

\begin{figure}[H]
		%\centering
		\includegraphics[scale=0.35]{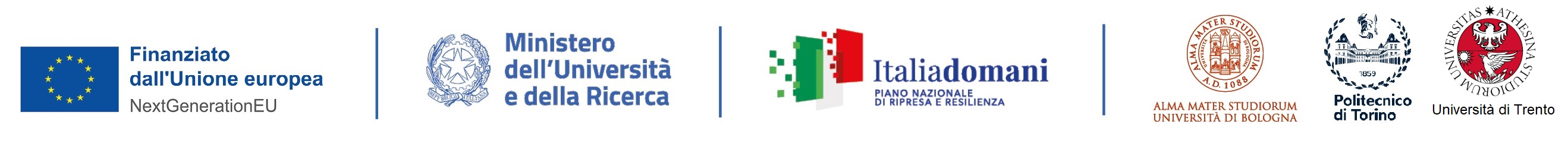}
			\end{figure}

All authors are members of INdAM-GNSAGA.

% \begin{figure}[H]
% 		%\centering
% 		\includegraphics[scale=0.35]{PRINBOTOTN.jpg}
% 			\end{figure}

% \bigskip

\end{document}